\documentclass[final]{siuro210301}

\usepackage{amsfonts}
\usepackage{graphicx}
\ifpdf
  \DeclareGraphicsExtensions{.eps,.pdf,.png,.jpg}
\else
  \DeclareGraphicsExtensions{.eps}
\fi
\graphicspath{{./images/}}  

\usepackage[T1]{fontenc}  
\usepackage{lmodern}  
\usepackage[utf8]{inputenc}  
\usepackage[english]{babel}  
\usepackage{microtype}
\usepackage{textcomp}  

\newcommand{\Reals}{\mathbb{R}}

\usepackage{listings}

\usepackage{enumitem}
\setlist[enumerate]{leftmargin=.5in}
\setlist[itemize]{leftmargin=.5in}

\headers{Efficient Interpolative Decompositions}{Rishi Advani and Sean O'Hagan}

\title{Efficient Algorithms for Constructing an Interpolative Decomposition}

\author{Rishi Advani\thanks{Cornell University (\email{ra534@cornell.edu}).}
\and Sean O'Hagan\thanks{University of Connecticut (\email{sean.ohagan@uconn.edu}).}}

\dedication{\small\textit{Project advisor: Akil Narayan\thanks{University of Utah.}}}

\ifpdf
\hypersetup{
  pdftitle={Efficient Algorithms for Constructing an Interpolative Decomposition},
  pdfauthor={Rishi Advani and Sean O'Hagan}
}
\fi

\begin{document}

\maketitle

\begin{abstract}
Low-rank approximations are essential in modern data science. The interpolative decomposition provides one such approximation. Its distinguishing feature is that it reuses columns from the original matrix. This enables it to preserve matrix properties such as sparsity and non-negativity. It also helps save space in memory. In this work, we introduce two optimized algorithms to construct an interpolative decomposition along with numerical evidence that they outperform the current state of the art.
\end{abstract}

\section{Introduction}

\subsection{Low-rank approximation}
As the dimensionality and size of modern datasets grow, numerical methods become increasingly important in making data analysis tractable. Low-rank approximations, in particular, are essential tools in numerical linear algebra and data science. They often make datasets simpler to work with, easier to understand, and more memory efficient. A comprehensive overview of low-rank approximation is given by \cite{halko_finding_2011}.

The Eckart–Young–Mirsky theorem shows that the truncated singular value decomposition (SVD) gives the closest rank-$k$ approximation to the original matrix in any unitarily invariant norm \cite{eckart_approximation_1936, mirsky_symmetric_1960}. One example of a unitarily invariant norm is the Frobenius norm, which we will be using in our analysis.

\subsection{Interpolative decomposition}
One downside of using the truncated SVD to form a low-rank approximation is that we are, in general, unable to preserve the original columns of the matrix. The one-sided interpolative decomposition (ID)\footnotemark, formally introduced to the literature by \cite{liberty_randomized_2007}, is able to do exactly that, and by doing so, is able to preserve more of the structure of the matrix (e.g., sparsity and non-negativity) By reusing columns of the original matrix, we are also able to save space in memory \cite{voronin_efficient_2017}.

\footnotetext{\cite{cheng_compression_2005} originally proposed a two-sided variant.}

\begin{definition}[Interpolative Decomposition]
Given an $m \times n$ matrix $A$, an $m \times k$ matrix $C$ whose columns constitute a subset of those of $A$, and a $k \times n$ matrix $Z$, such that
\begin{itemize}
    \item some size-$k$ subset of the columns of $Z$ form the $k \times k$ identity matrix, and
    \item no entry of $Z$ has absolute value greater than $2$,
\end{itemize}
$CZ$ is an \emph{interpolative decomposition} of $A$.
\end{definition}

Note that our definition is equivalent to the ``weaker'' form given by \cite{liberty_randomized_2007}, but without specific guarantees on the accuracy of the approximation.
Further variations on the ID have been proposed \cite{ari_probabilistic_2012}, and several optimized algorithms for computing IDs have been designed \cite{voronin_efficient_2017, malik_fast_2019, lucas_parallel_2014, woolfe_fast_2008, martinsson_randomized_2011}.
The ID has seen applications in computational electromagnetic problems, in particular to ``method of moments'' systems \cite{pan_fast_2012, si-lu_huang_efficient_2016, pan_preconditioning_2013}.

\cite{martinsson_id_2014} have developed a software package containing Fortran implementations of two algorithms (one deterministic and one randomized) for computing IDs. The implementations are based on work done by various authors \cite{cheng_compression_2005, liberty_randomized_2007, woolfe_fast_2008, martinsson_randomized_2011}. A Python wrapper of this package is included in the SciPy \cite{scipy_10_contributors_scipy_2020} library. To the best of the authors' knowledge, these algorithms represent the current state of the art for computing a one-sided ID. From here on, we will refer to the deterministic implementation as SciPy ID and the randomized one as SciPy RID.

\subsection{Contributions}
We contribute two algorithms: Optim ID and Optim RID. The first is a Python implementation of the algorithm behind SciPy ID that proves to be more efficient. The second is the result of applying the general approach of Optim ID to the method of \cite{advani_random_2020} to create a state-of-the-art randomized method for computing an ID.

\subsection{Outline}
\Cref{algorithm} contains descriptions and analyses of our algorithms.
\Cref{numerical_results} contains numerical evidence that our algorithms successfully construct IDs and outperform the methods provided by the SciPy library.
\Cref{conclusion} contains concluding remarks and potential directions for future research.

\section{Algorithms} \label{algorithm}
Let $A \in \Reals^{m \times n}$ with $m \leq n$ be a matrix\footnote{If the matrix $A$ has $m > n$, we apply our method to $A^T$ and transpose the resulting decomposition to get the ``dual'' of an ID (i.e., the product of a matrix consisting of rows of $A$ and a matrix containing the identity with small entries).}
of rank at least $k$. We present algorithms to construct a rank-$k$ ID of $A$ below.

\subsection{Deterministic algorithm}
We describe the algorithm behind SciPy ID and provide our implementation, Optim ID. At a high level, the algorithm uses column-pivoted QR to select columns for $C$ and then computes $Z$ via least-squares. Our implementation is displayed as Algorithm~\ref{alg_optim_id}.

First, we compute the column-pivoted QR factorization, $AP = QR$, where $Q$ is orthogonal, $R$ is upper triangular, and $P$ is a permutation matrix. Let $Q_k$ denote the submatrix of $Q$ consisting of the first $k$ columns of $Q$. Let $R_k$ denote the submatrix of $R$ consisting of only the entries in both the first $k$ columns and rows of $R$. Let $P_k$ denote the first $k$ columns of $P$. Let $C=AP_k$. As long as the rank of the matrix $A$ is at least $k$, then $C$ is full rank.

To produce an accurate approximation, we aim to find the matrix $Z$ that minimizes the following error: $\lVert A - CZ \rVert_F$. This is a least-squares problem.

The solution to a least-squares problem is given by the exact solution $Z$ to the normal equations,
\begin{equation}
    C^{T}CZ = C^{T}A \,. \label{normal_eq}
\end{equation}
By the properties of the QR factorization, we have
\[C = AP_k = Q_k R_k \,.\]
We can then simplify \cref{normal_eq}:
\[R_k^T Q_k^T Q_k R_k Z = C^T A\]
Since $Q_k$ has orthonormal columns, we can further simplify:
\begin{equation}
    R_k^T R_k Z = C^T A \label{final_eq}
\end{equation}

The matrix $C$ is full rank, so the diagonal of $R_k$ has no zeros. Furthermore, $R_k$ is triangular, so it is nonsingular. Finally, we note that $R_k^T R_k$ is positive definite.

We then use back-/forward-substitution on the triangular linear systems to efficiently solve \cref{final_eq} for $Z$. The approximation is given by $A \approx CZ$.

\begin{lstlisting}[caption=Optim ID, label=alg_optim_id, float=tphb]
def optim_id(A, k):
    _, R, P = scipy.linalg.qr(A, 
            pivoting=True, 
            mode='economic', 
            check_finite=False)
    R_k = R[:k,:k]
    cols = P[:k]
    C = A[:,cols]
    Z = scipy.linalg.solve(R_k.T @ R_k, 
            C.T @ A, 
            overwrite_a=True, 
            overwrite_b=True, 
            assume_a='pos')
    approx = C @ Z
    return approx, cols, Z
\end{lstlisting}

\subsection{Randomized algorithm}
Here, we take the column sampling idea introduced by \cite{advani_random_2020} and apply to it the ideas behind Optim ID to obtain the algorithm Optim RID. At a high level, the algorithm randomly samples $p$ columns from $A$, uses column-pivoted QR to select $k$ of those $p$ columns for $C$, and then computes $Z$ via least-squares. Our implementation is displayed as Algorithm~\ref{alg_optim_rid}.

If we were to naively sample random columns of $A$, we would likely not capture the full range of the matrix. To help ensure that we do, we oversample to a certain degree, depending on the desired rank of the approximation. This greatly increases the probability that the sampled vectors span a large portion of the range of $A$. By default, we take our oversampling parameter to be $p = 1.2k$, but this can be adjusted as necessary. We randomly sample (without replacement) $p$ columns from $A$. Let the matrix formed by these columns be denoted by $A_S$.

We compute the column-pivoted QR factorization, $A_S P = QR$, where $Q$ is orthogonal, $R$ is upper triangular, and $P$ is a permutation matrix. Let $Q_k$ denote the submatrix of $Q$ consisting of the first $k$ columns of $Q$. Let $R_k$ denote the submatrix of $R$ consisting of only the entries in both the first $k$ columns and rows of $R$. Let $P_k$ denote the first $k$ columns of $P$. Let $C=A_S P_k$.

To produce an accurate approximation, we aim to find the matrix $Z$ that minimizes the following error: $\lVert A - CZ \rVert_F$.

As before, the matrix $Z$ is given by the solution to \cref{final_eq}:
\[R_k^T R_k Z = C^T A\]

Since $C$ may be rank deficient, we are unable to assert that $R_k^T R_k$ is positive definite. Instead, we fall back on the weaker property that $R_k^T R_k$ is a symmetric matrix. We then use the diagonal pivoting method \cite[routine \texttt{dsysv}]{dodge_lapack_1992} to solve \cref{final_eq}. The approximation is given by $A \approx CZ$.

\begin{lstlisting}[caption=Optim RID, label=alg_optim_rid, float=tphb]
rng = numpy.random.default_rng()
def optim_rid(A, k):
    oversampling = int(0.2 * k)
    p = k + oversampling
    idx = rng.choice(A.shape[1], 
            replace=False, 
            size=p)
    AS = A[:,idx]
    _, R, P = scipy.linalg.qr(AS, 
            pivoting=True, 
            mode='economic', 
            check_finite=False)
    R_k = R[:k,:k]
    _cols = P[:k]
    cols = idx[_cols]
    C = AS[:,_cols]
    Z = scipy.linalg.solve(R_k.T @ R_k, 
            C.T @ A, 
            overwrite_a=True, 
            overwrite_b=True, 
            assume_a='sym')
    approx = C @ Z
    return approx, cols, Z
\end{lstlisting}

\subsubsection{Comparison with SciPy RID}
The key to both algorithms is finding a smaller matrix with which to compute the ID. SciPy RID relies on a composition of a random transform, a fast Fourier transform, and column sampling. Optim RID uses column sampling directly, which we will see from the results gives us increased efficiency at the cost of reduced accuracy on very sparse datasets.

\section{Numerical results} \label{numerical_results}

\subsection{ID Properties}
To verify that the decompositions given above are indeed IDs, we need to check that $C$ is comprised solely by columns of $A$ and that entries of $Z$ are bounded by 2. Mathematically, the first claim holds by construction.

The second claim is shown to hold in practice through our numerical experiments. Optim ID was able to bound the entries of $Z$ by 2 on all datasets. Optim RID was able to bound the entries on all dense datasets, but not on sparse ones. The maximum entries of $Z$ for a fixed rank of $k=190$ are displayed in \cref{table_max_entry}. For the randomized algorithms, the number displayed is the mean over 10 iterations of the algorithm. Note in particular that Optim RID performs poorly on the sparsest dataset, Sparse1. Full results can be found in \cref{full_results}.

\begin{table}[hpt]
\footnotesize\caption{Max Entries of $Z$ for Rank-190 Approximations}\label{table_max_entry}
\centering
\begin{tabular}{|c||c|c|c|c|c|}
    \hline
    Dataset & SciPy ID & Optim ID & SciPy RID & Optim RID \\
    \hline
    Boolean & 1 & 1 & 1 & 1 \\
    Gaussian & 1 & 1 & 1 & 1 \\
    Uniform & 1 & 1 & 1 & 1 \\
    MNIST & 1 & 1 & 1.010 & 1.004 \\
    Fashion & 1 & 1 & 1.027 & 1.022 \\
    Sparse1 & 1 & 1 & 1.058 & 167.33 \\
    Sparse2 & 1 & 1 & 1.011 & 2.39 \\
    Sparse3 & 1 & 1 & 1.008 & 3.32 \\
    \hline
\end{tabular}
\end{table}

\subsection{Performance}
We tested the accuracy and computational efficiency of Optim ID and Optim RID against that of SciPy ID and SciPy RID. We used the SVD as a baseline in all tests. In order to ensure robustness, we tested on a variety of datasets; the precise details of each dataset can be found in \cref{full_results}. Most of the datasets are roughly $1000 \times 1000$ in matrix form. MNIST and Fashion-MNIST are somewhat larger and Sparse3 is somewhat smaller.
In each test, we measured the relative error (with respect to the Frobenius norm), execution time, and ability to bound the entries of the matrix $Z$.

In \cref{table_error,table_time}, we list the error and time results for a fixed rank of $k=190$. In each row, the entry for the algorithm with the best result is shown in bold. For the randomized algorithms, the number displayed is the mean over 10 iterations of the algorithm.
Note that in \cref{table_error}, SciPy ID and Optim ID share a column, as both algorithms construct the same decomposition and hence have the same error. Also, the SVD is not marked in bold, as it necessarily represents the optimal approximation by the Eckart–Young–Mirsky theorem; it is only shown as a baseline.

\begin{figure}[htpb]
\centering
\includegraphics[width=\linewidth]{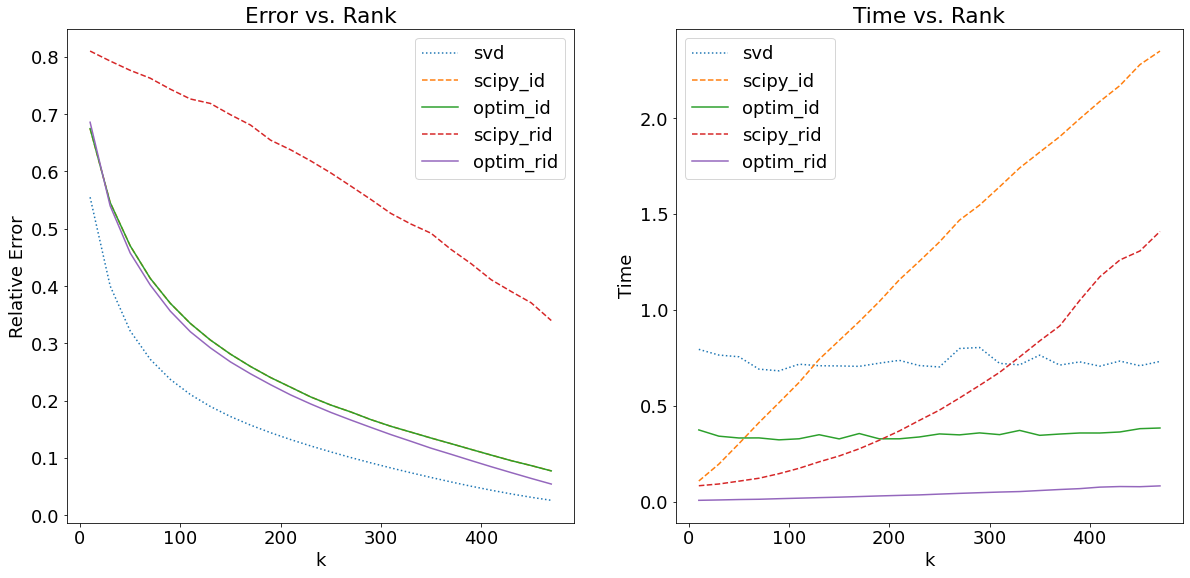}
\caption{Relative Error and Time (s) on MNIST. Note that the SciPy ID line coincides with the Optim ID line in the left graph.}
\label{mnist_mini}
\end{figure}

In \cref{mnist_mini}, we show the performance of the algorithms on the MNIST dataset for various ranks $k$.
Full results can be found in \cref{full_results}.

\begin{table}[hpt]
\footnotesize\caption{Relative Error for Rank-190 Approximations}\label{table_error}
\centering
\begin{tabular}{|c||c|c|c|c|c|}
    \hline
    Dataset & SVD & SciPy/Optim ID & SciPy RID & Optim RID \\
    \hline
    Boolean & .467 & \textbf{.553} & 1.401 & .554 \\
    Gaussian & .660 & \textbf{.776} & 1.988 & .782 \\
    Uniform & .331 & \textbf{.390} & 1.003 & .392 \\
    MNIST & .144 & .240 & .664 & \textbf{.228} \\
    Fashion & .140 & .215 & .613 & \textbf{.200} \\
    Sparse1 & .020 & \textbf{.022} & .086 & .783 \\
    Sparse2 & .481 & \textbf{.540} & 1.767 & .697 \\
    Sparse3 & .278 & \textbf{.320} & 1.167 & .535 \\
    \hline
\end{tabular}
\end{table}

\begin{table}[hpt]
\footnotesize\caption{Execution Time (s) for Rank-190 Approximations}\label{table_time}
\centering
\begin{tabular}{|c||c|c|c|c|c|c|}
    \hline
    Dataset & SVD & SciPy ID & Optim ID & SciPy RID & Optim RID \\
    \hline
    Boolean & .094 & .189 & .039 & .063 & \textbf{.011} \\
    Gaussian & .088 & .181 & .036 & .061 & \textbf{.010} \\
    Uniform & .088 & .195 & .038 & .060 & \textbf{.010} \\
    MNIST & .721 & 1.042 & .327 & .654 & \textbf{.228} \\
    Fashion & .707 & 1.036 & .335 & .613 & \textbf{.200} \\
    Sparse1 & .177 & .325 & .071 & .081 & \textbf{.015} \\
    Sparse2 & .088 & .167 & .036 & .051 & \textbf{.010} \\
    Sparse3 & .030 & .056 & .014 & .027 & \textbf{.006} \\
    \hline
\end{tabular}
\end{table}

\section{Conclusion} \label{conclusion}
We have introduced two new algorithms for computing low-rank approximations, Optim ID and Optim RID. The algorithms SciPy ID and Optim ID both compute the same approximation, but in almost all cases, Optim ID computes the approximation in less time. It is unclear why this is the case, as SciPy ID and Optim ID both rely on the same underlying algorithm. The main implementation difference is that SciPy ID is written in Fortran and Optim ID is written in Python.

Optim RID is always the most efficient method, and it is more accurate than SciPy RID on all but one dataset. It often has accuracy near that of the deterministic algorithms and occasionally even surpasses it. Specifically, Optim RID outperforms the deterministic methods on the two real-world dense datasets we tested on, MNIST and Fashion-MNIST.

The algorithms presented in this paper consistently outperform the methods provided in the SciPy library, and in doing so, achieve a new state of the art. We hope these novel methods make analysis of large datasets more tractable.

\subsection{Future work}
We list several possible directions for future research.

For our experiments, we used an oversampling parameter of $0.2k$ (i.e., we sampled $1.2k$ columns when we wanted an approximation of rank $k$). This value proved sufficient for most datasets, but our randomized algorithm was unable to sufficiently bound the entries of the matrix $Z$ when tested on sparse datasets. We suspect that the algorithm would be able to scale to more sparse datasets if the oversampling parameter was partially determined by the sparsity of the dataset.

One further optimization that we could have used, but did not have the time to properly implement and test, would be halting the computation of the QR factorization after obtaining a set of $k$ linearly independent vectors. In our current implementation, we compute the full factorization, then throw out the last $n-k$ columns.

While we tested on dense and sparse datasets, we did not experiment to find the ``critical point'' of sparsity (if such a value exists) where the accuracy of SciPy RID overtakes that of Optim RID. This would be useful in constructing a composite algorithm that executes a particular sub-algorithm depending on the detected level of sparsity of the dataset.

As our algorithms are optimized for dense matrices, it would be interesting to see if an improved algorithm for sparse matrices could be designed. One promising idea is to use an iterative implementation of the conjugate gradient method. On very large, sparse matrices, this may result in improved performance.

Finally, it is possible to define the ID for matrices with complex entries, and much of the theoretical analysis does not change. We did not test our methods on datasets with complex entries, and some modifications may be necessary for the code to compile, but it would be interesting to see how the methods perform on such data.

\appendix
\section{Full Results} \label{full_results}

The algorithms were tested using the following ranks $k$:
\[\lbrace \, 10, 30, 50, \dots, 450, 470 \, \rbrace\]

The following dense datasets were tested on:
\begin{itemize}
    \item Boolean:
    a matrix of dimensions $784 \times 1000$ with entries randomly sampled from the set $\lbrace \, 0,1 \, \rbrace$
    
    \item Gaussian:
    a matrix of dimensions $784 \times 1000$ with entries sampled from a standard normal distribution
    
    \item Uniform:
    a matrix of dimensions $784 \times 1000$ with entries sampled from a uniform distribution over the half-open interval $[0,1)$
    
    \item MNIST:
    the first 5000 images in the training data of MNIST \cite{lecun_gradient-based_1998} (each image is flattened to a vector)
    
    \item Fashion:
    the first 5000 images in the training data of Fashion-MNIST \cite{xiao_fashion-mnist_2017} (each image is flattened to a vector)
\end{itemize}

The sparse datasets tested on are all from the SuiteSparse Matrix Collection \cite{davis_university_2011}. Dataset-specific information can be found in \cref{sparse_table}.

\begin{table}[htbp]
\footnotesize\caption{Sparse Dataset Information}\label{sparse_table}
\centering
\begin{tabular}{|c||c|c|c|c|c|c|c|} 
\hline
Dataset & ID & Name & Rows & Cols & Nonzeros & Type & Sparsity \\
\hline
Sparse1 & 1 & \verb|1138_bus| & 1138 & 1138 & 4054 & Power Network & .00313 \\
Sparse2 & 2888 & \verb|Vehicle_10NN| & 846 & 846 & 10894 & Weighted Graph & .01522 \\
Sparse3\footnotemark & 2885 & \verb|Spectro_NN| & 531 & 531 & 7422 & Weighted Graph & .02632 \\
\hline
\end{tabular}
\end{table}

\footnotetext{The Sparse3 dataset was tested only on ranks $k$ less than 400 because of its relatively small size.}

Our algorithms, Optim ID and Optim RID, are implemented in Python (v3.8.5) and use various methods from the SciPy (v1.5.2) and NumPy \cite{harris_array_2020} (v1.19.2) libraries (e.g., the SciPy implementation of QR decomposition). For our baseline method, we used the NumPy implementation of the SVD.

All tests were run on the same machine. The specifications are provided below:
\begin{itemize}
    \item Processor: Intel(R) Core(TM) i7-10700K CPU @ 3.80 GHz, 8 cores, 16 threads
    \item RAM: 16.00 GB
    \item Cache: 16MB Intel Smart Cache
    \item OS: 64-bit Windows 10 Pro
\end{itemize}

The following graphs represent the results for the listed ranks on each dataset.

\begin{figure}[htpb]
\centering
\includegraphics[width=\linewidth]{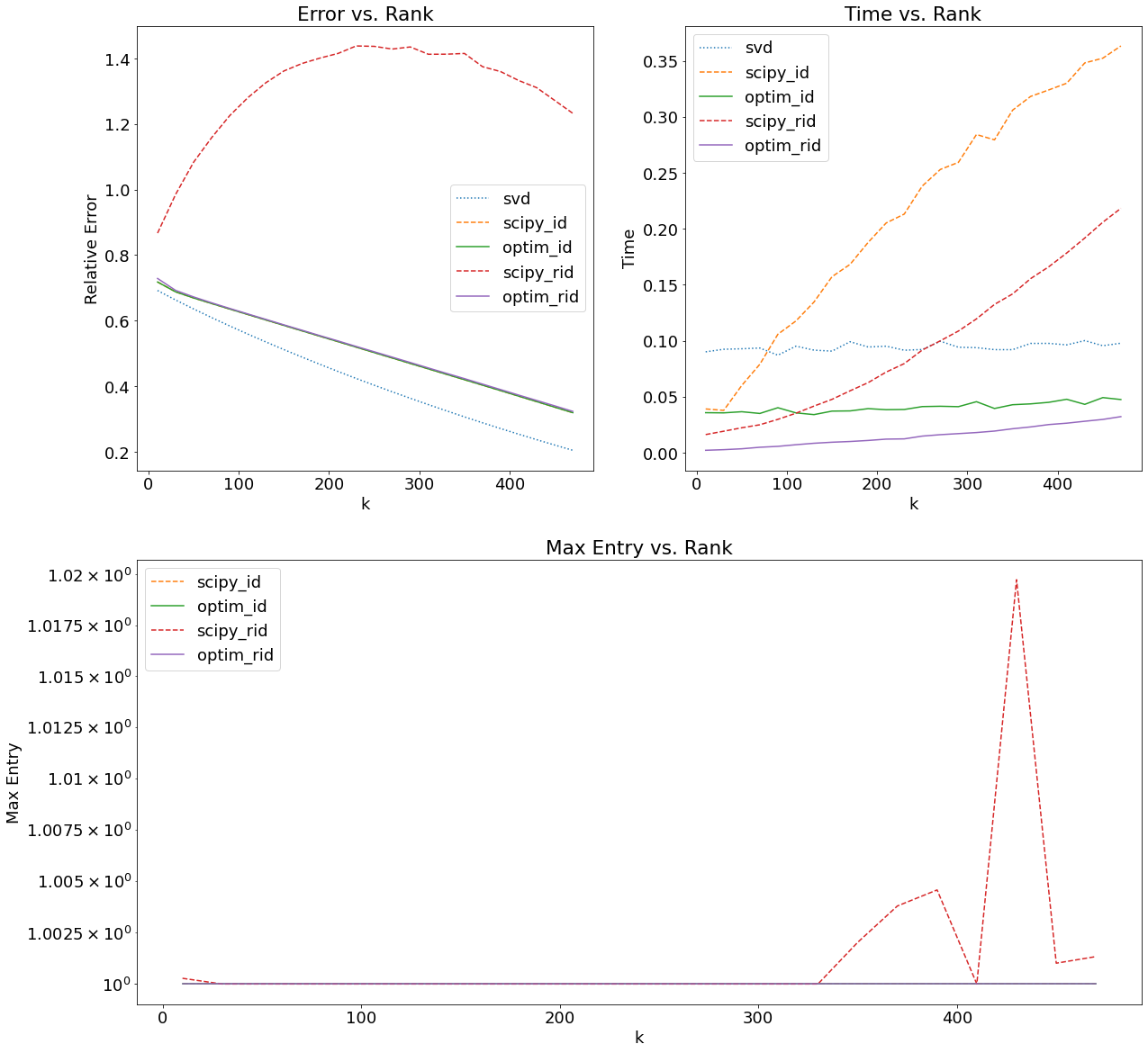}
\caption{Boolean}
\end{figure}

\begin{figure}[htpb]
\centering
\includegraphics[width=\linewidth]{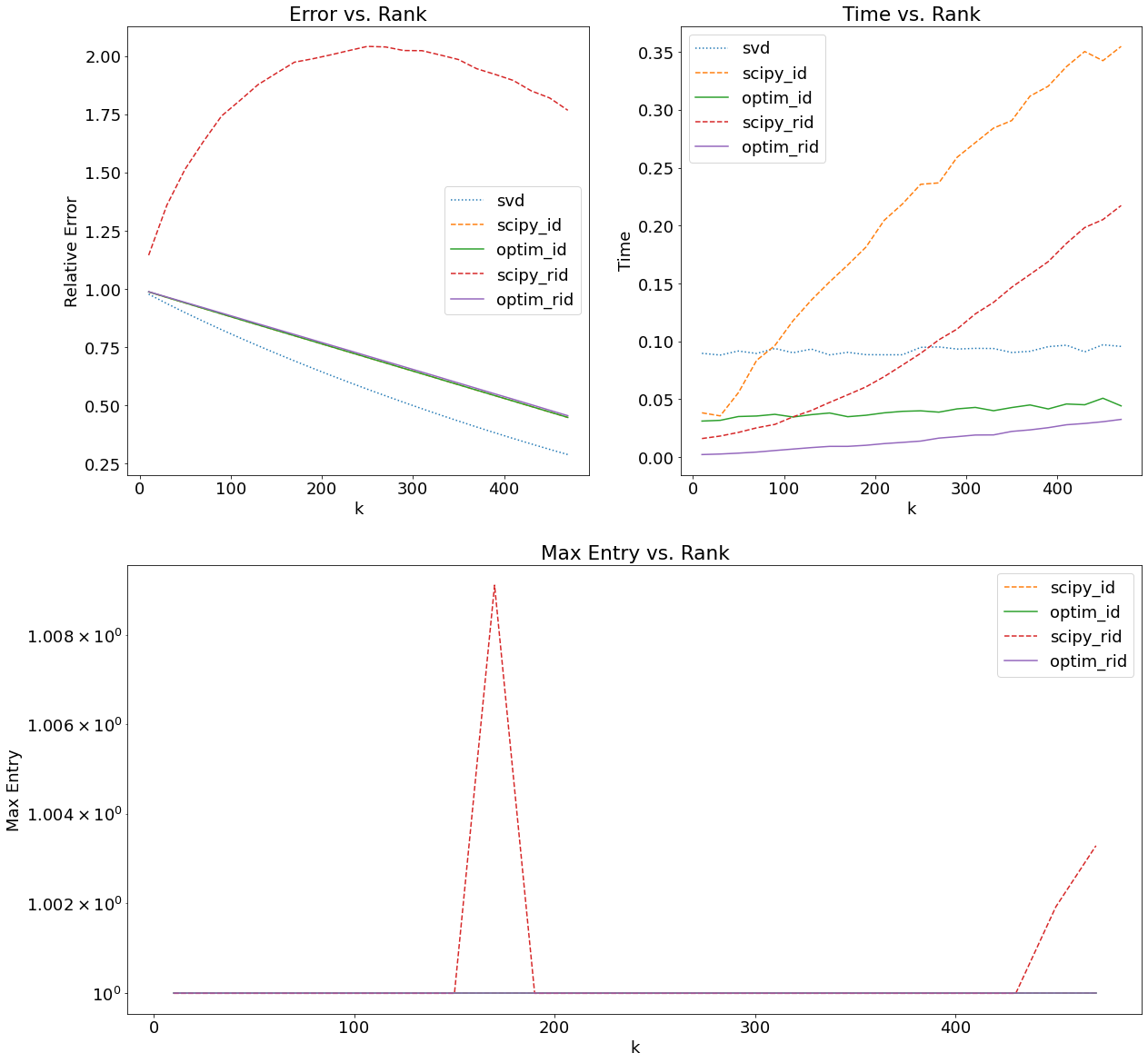}
\caption{Gaussian}
\end{figure}

\begin{figure}[htpb]
\centering
\includegraphics[width=\linewidth]{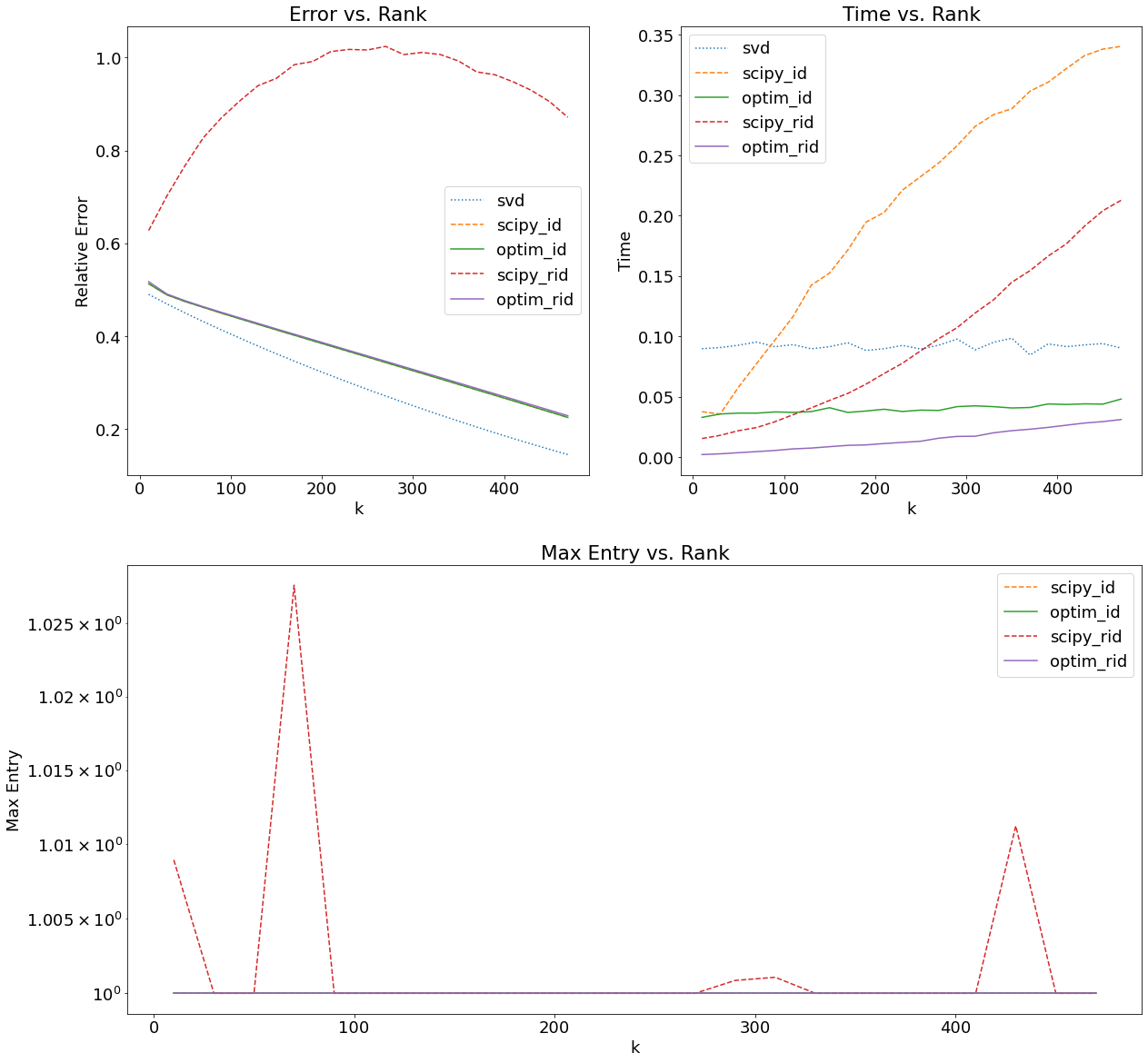}
\caption{Uniform}
\end{figure}

\begin{figure}[htpb]
\centering
\includegraphics[width=\linewidth]{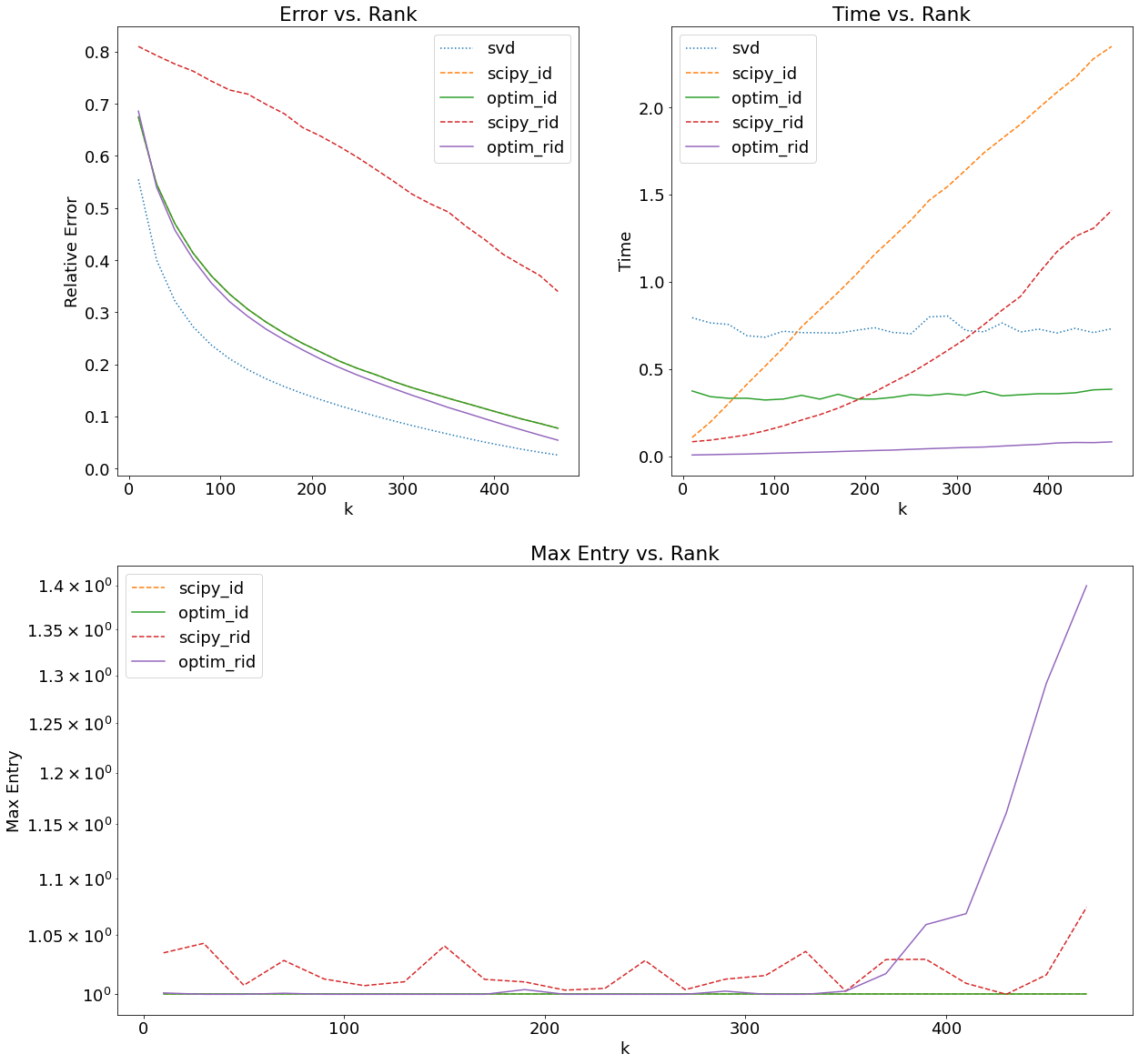}
\caption{MNIST}
\end{figure}

\begin{figure}[htpb]
\centering
\includegraphics[width=\linewidth]{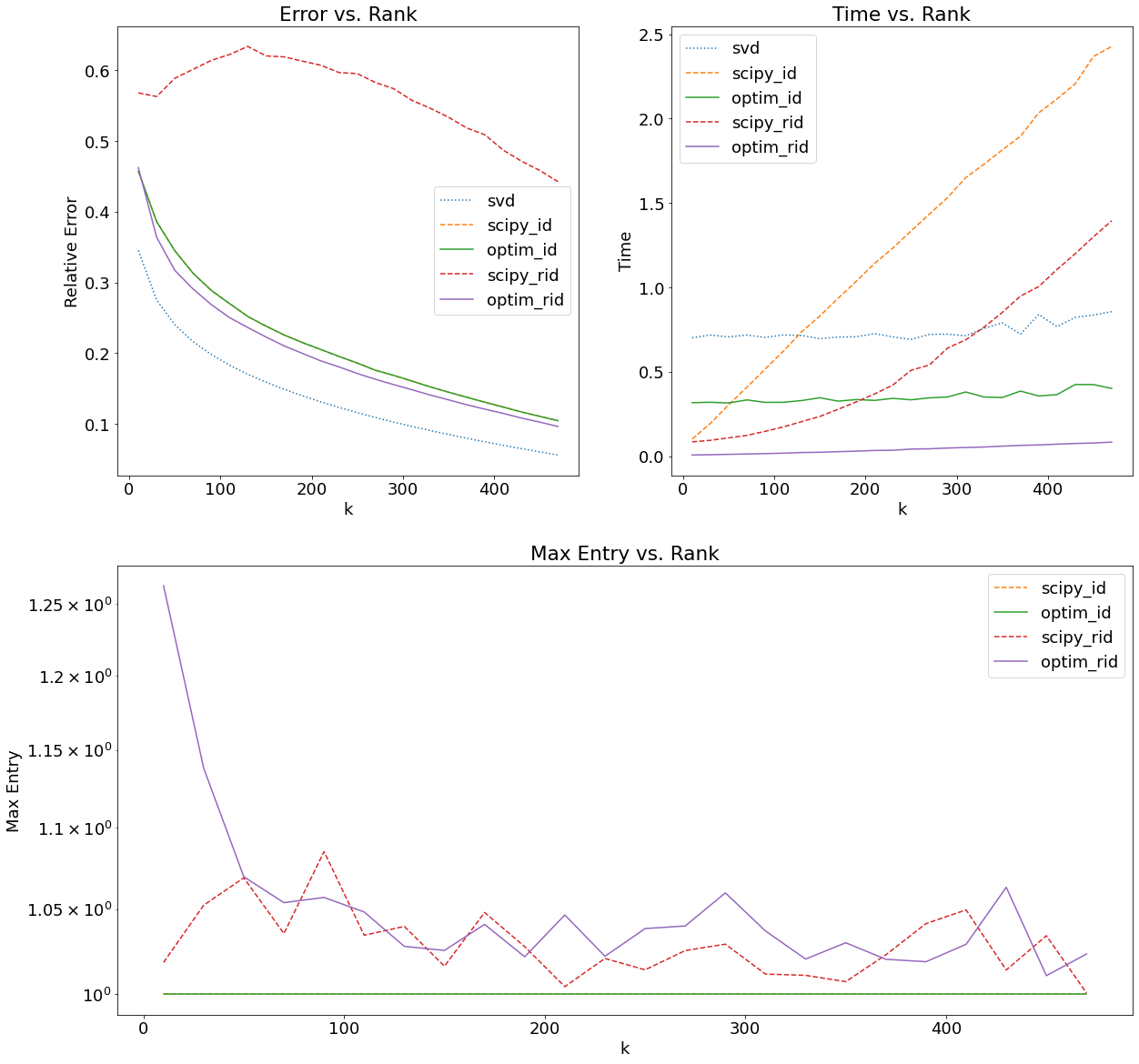}
\caption{Fashion}
\end{figure}

\begin{figure}[htpb]
\centering
\includegraphics[width=\linewidth]{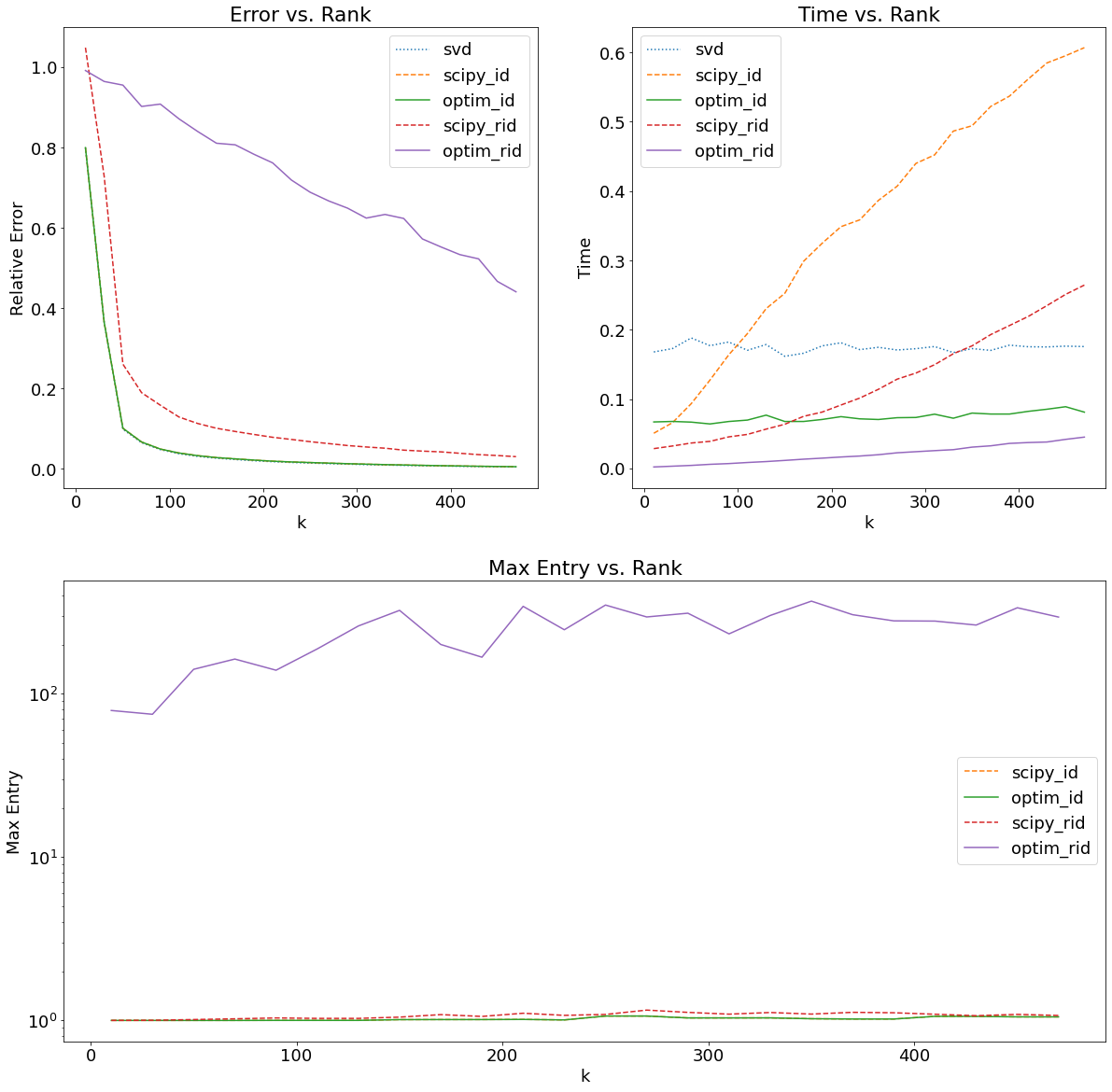}
\caption{Sparse1}
\end{figure}

\begin{figure}[htpb]
\centering
\includegraphics[width=\linewidth]{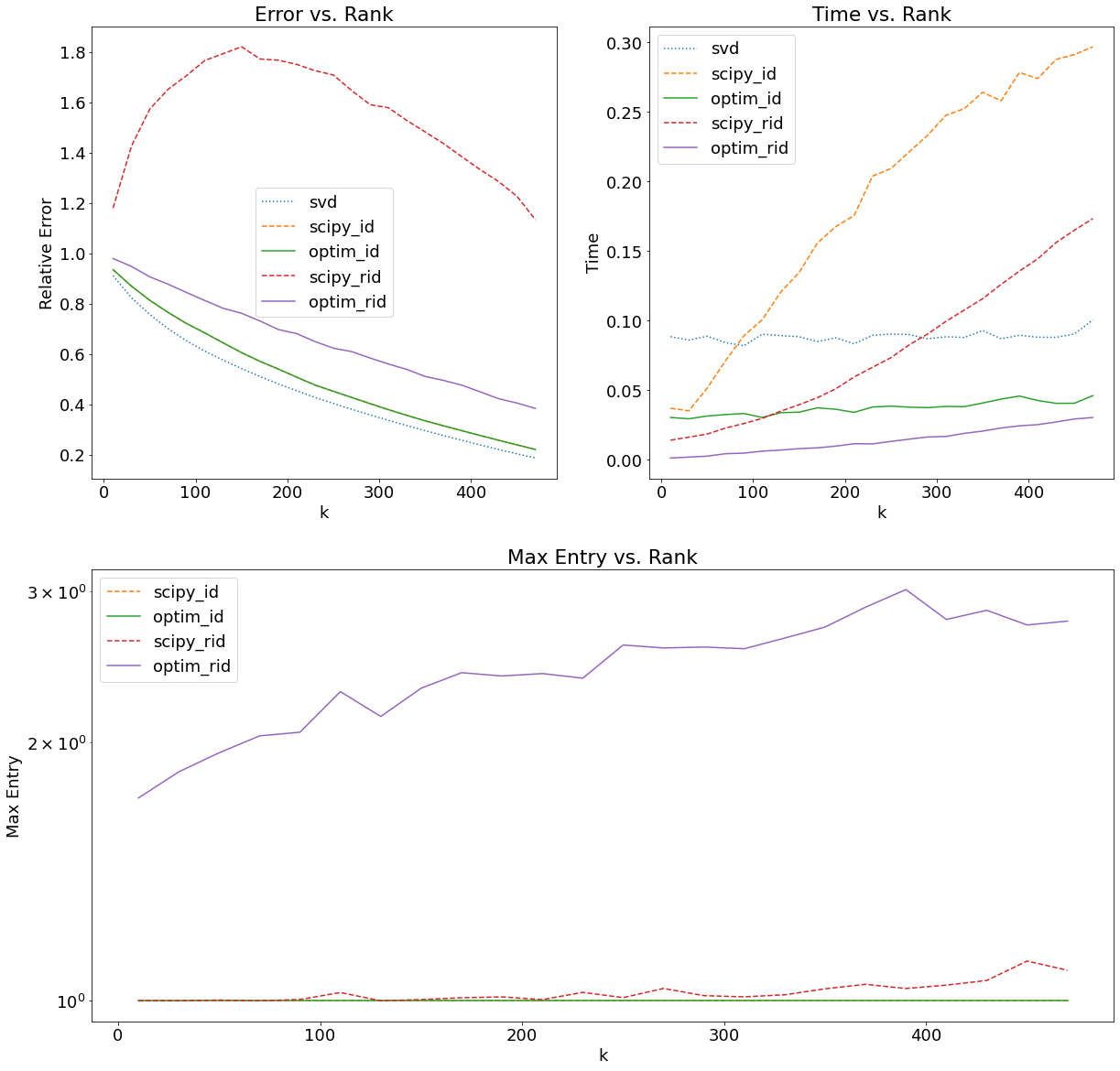}
\caption{Sparse2}
\end{figure}

\begin{figure}[htpb]
\centering
\includegraphics[width=\linewidth]{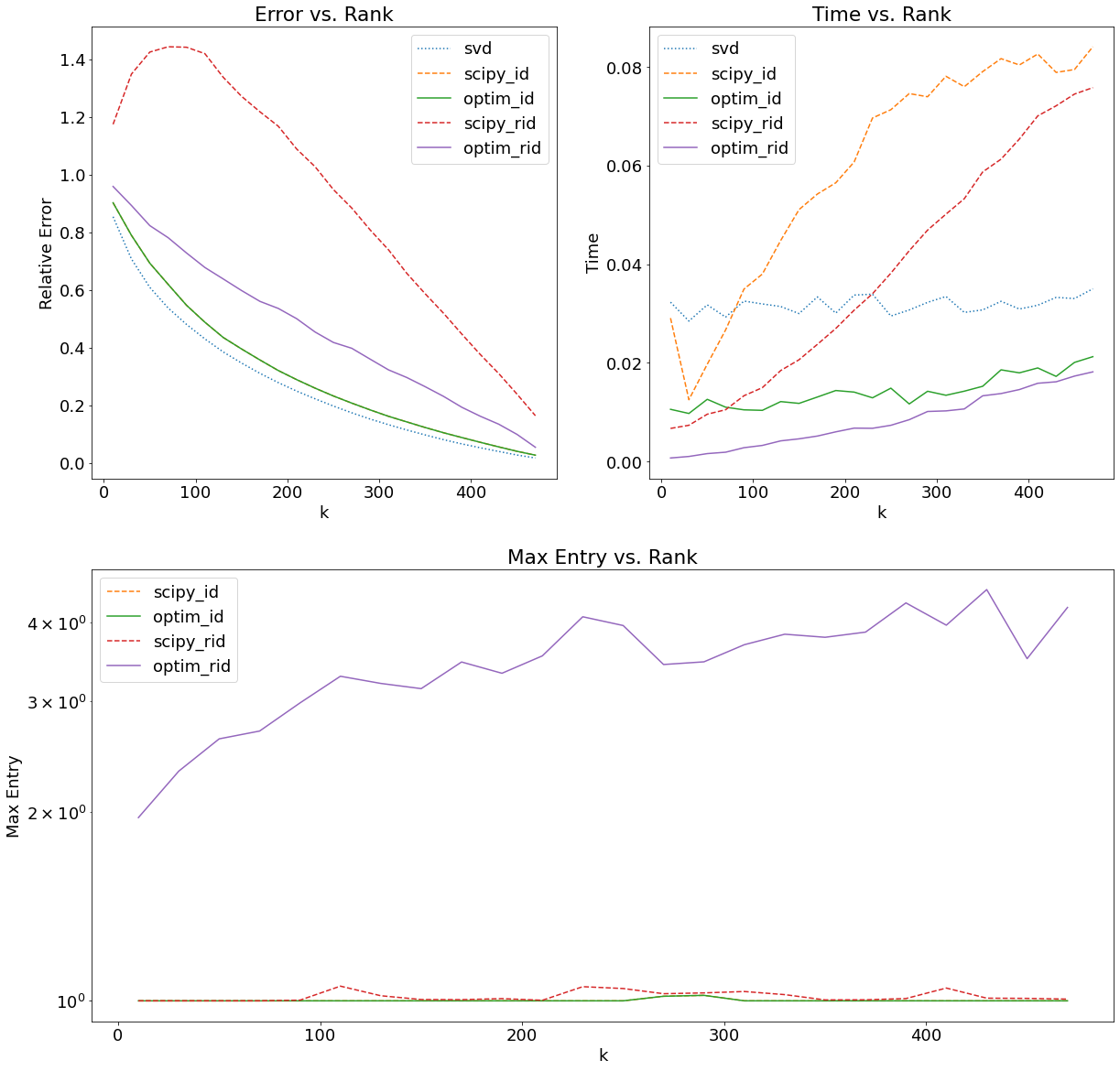}
\caption{Sparse3}
\end{figure}

\section*{Acknowledgments}
This paper is a continuation of work done with Madison Crim as undergraduate research fellows at ICERM. We would like to thank the organizers of Summer@ICERM 2020; our primary advisor, Akil Narayan; our secondary advisor, Yanlai Chen; and our TAs, Justin Baker and Liu Yang. We would also like to thank Karen Zhou for proofreading the final document.

\bibliographystyle{siamplain}
\bibliography{main}

\end{document}